\documentclass[reqno]{article}

\usepackage{fullpage}

\usepackage{amssymb,amsmath,epsfig,graphics,color,float,subfigure}
\usepackage[ansinew]{inputenc}
\newtheorem{theorem}{Theorem}[section]
\newtheorem{proposition}[theorem]{Proposition}
\newtheorem{lemma}[theorem]{Lemma}

\newtheorem{remark}[theorem]{Remark}

\def\ds{\displaystyle}

\title{\bf Extensions of the chemostat model with flocculation}
\author{{\bf R. Fekih-Salem$^{a,b}$, J. Harmand$^{b,c}$,
    C. Lobry$^{b}$, A. Rapaport$^{a,b}$\thanks{corresponding
      author. Tel: +33.4.99.61.26.52 Fax: +33.4.67.52.14.27. E-mails:
      radhouene.fs@gmail.com (R. Fekih-Salem), harmand@supagro.inra.fr
      (J. Harmand), claude.lobry@inria.fr (C. Lobry), 
      rapaport@supagro.inra.fr (A. Rapaport), tewfik.sari@irstea.fr
      (T. Sari)} \, and T. Sari$^{b,d}$}\\
\small{$^{a}$ UMR INRA-SupAgro MISTEA, 1 p. Viala, 34060 Montpellier,
  France}\\[-1mm]
\small{$^{b}$ EPI INRA-INRIA MODEMIC, route des Lucioles, 06902 Sophia-Antipolis, France}\\[-1mm]
\small{$^{c}$ INRA LBE, Avenue de Etangs, 11100 Narbonne, France}\\[-1mm]
\small{$^{d}$ Irstea, UMR ITAP, 361 rue Jean-François Breton, 34196 Montpellier}}

\date{\today}

\begin{document}
\maketitle
\begin{abstract}
In this work, we study a model of the chemostat where the species are
present in two forms, isolated and aggregated individuals, such as
attached bacteria or bacteria in flocks. We show that our general
model contains a lot of models that were previously considered in the
literature. Assuming that flocculation and deflocculation dynamics are
fast compared to the growth of the species, we construct a reduced chemostat-like model in which both the growth functions and the apparent dilution rate depend on the density of the species. We also show that such a model involving monotonic growth rates may exhibit bi-stability, while it
may occur in the classical chemostat model, but when the growth rate is non monotonic.\\

\end{abstract}
\section{Introduction}

In culture of micro-organisms, the attachment of microbial individuals occurs frequently. The attachment can be either a ``wall attachment'' such as in the growth of biofilms or simply an aggregation such as in the formation of flocks or granules \cite{C95,TJF99}.
Flock or granule formation has a direct impact
on growth dynamics, as the access to the substrate is limited for micro-organisms inside such structures.
The mechanisms of attachment and detachment result from the coupling of hydrodynamics conditions and biological properties, but are not yet completely understood at the level of microbial individuals. Several attempts of computer models, using individual based representations, have been proposed and are under investigation for the simulation of these phenomenons, cf. for instance \cite{iwa06,SS06}.
At a macroscopic level, substrate limitation can be measured experimentally in biofilms or flocks \cite{CCLWL07,LiBi04,TH71}.
A rough representation, suited to the
macroscopic level, consists in splitting the overall biomass into two parts: a ``planktonic biomass'', composed of free individuals and an ``attached biomass'' composed of individuals that are tied together \cite{KS85}.
This consideration leads to a significant change on the performances predicted by the models, compared to purely planktonic cultures. In a chemostat-like device,
planktonic cells are expected to consume easily the substrates necessary for their growth, but are
more keen to be carried out by the flow. On the contrary, cells
among aggregates or biofilms have a more difficult access to the
resources of the bulk fluid, but are more resistant to detachment induced by the hydrodynamical conditions.
Therefore, mathematical models are expected to understand and predict the issues of these trade-offs.
Several extensions of the well-known chemostat model \cite{chem}, considering two compartments of free and attached biomass for each species have been proposed and studied in the literature.
In models with wall attachment, attached biomass is assumed to be fixed while detached individuals return directly to the planktonic compartment \cite{PW99,SS00}.
In models with aggregation, aggregates are carried out by the flow but bacteria inside flocks are assumed to have no or reduced access to bulk resources \cite{HLH07,HR08}.

Literature reports flocculation time scales of
the order of 1 to 10 min \cite{ding,walberg} to be compared with bacterial
growth times of 1 h to 1 day, and with retention times of
a few hours to a few days.
Thus, considering that attachment and detachment processes may be fast compared to biological time, it is shown in \cite{HR08}
that the reduced dynamics of such systems amounts to have a single biomass compartment for each strain but with a density dependent growth rate.
This justifies the consideration of density dependent growth functions in the chemostat model, as already introduced in the literature in the field of mathematical ecology \cite{AG89}
 or waste-water process engineering \cite{HG}.
In \cite{LAS06,LH06,LMR04,LMR06}, it has been shown that this could
lead to the coexistence of several species in competition on a same
limiting resource, thus invalidating the Competitive Exclusion
Principle \cite{H60} (different mechanisms with considerations on
the nutrient uptake could also lead to species coexistence \cite{APW03}).

In \cite{HR08}, the aggregates are assumed to have no biological growth (i.e. the attachment process is the only source of increase of the attached biomass). Aggregates are also assumed to be washed-out with the same dilution rate than planktonic cells.
On the opposite, in wall attachment models, the attached biomass is not washed out at all.
We believe that these two opposite cases (same dilution rate than planktonic biomass or no dilution rate) are too extreme to be fully realistic. In this paper, we revisit the chemostat model with two compartments, planktonic and aggregated biomass, but assuming that each biomass has its own growth rate and apparent dilution rate.
This generalizes the two kind of models that we mentioned previously.

The paper is organized as follows. In the Section 2, the general model of the system under interest is presented. In particular, for specific choices of different kinetics and mortality terms, it is shown that this model captures in fact many models of literature. In Section 3, this general model is reduced assuming the attachment/detachment processes are fast with respect to others. In the Section 4, the reduced model is analyzed for a specific class of models while its extension to the multi-species case and its analysis are carried out in Sections 5 and 6, respectively. Finally, discussion and conclusions are drawn in the last section.

\section{Modeling flocks or aggregates in the chemostat}
Consider the following model
of the chemostat in which a population of microorganisms compete
for a single growth-limiting substrate \cite{chem,monod}:
$$\left\{
\begin{array}{lcl}
\dot{S} & = &  D(S_{in}-S)-k{\mu(S)}x\\
\dot{x} & = & (\mu(S)-D_x)x
\end{array}
\right..
$$
In these equations, $S(t)$ denotes the concentration of the substrate at time $t$;
$x(t)$ denotes the concentration of the  population of microorganisms at time $t$;
$\mu(\cdot)$ represents the per-capita growth rate of the population and so $Y=1/k$ is
the growth yield;
$S_{in}$ and $D$ denote, respectively, the concentration of substrate in the feed
bottle and the dilution rate of the chemostat;
$D_x$ represents the removal rate of the population.

{Assume that the species is present in two forms:}
isolated or planktonic bacteria, of density $u$,
and attached bacteria or flocks of bacteria, of density $v$.
Isolated bacteria and flocks can stick together
to form new flocks,  with rate $\alpha(\cdot)u$,
and flocks can split and liberate isolated bacteria,
with rate $\beta(\cdot)v$:
$$u\xrightarrow{\alpha(\cdot)u}v,\qquad u\xleftarrow{\beta(\cdot)v}v.$$
One obtains the following equations :
\begin{equation}\label{genmod}
\left\{
\begin{array}{lcl}
\dot{S}&=&D(S_{in}-S)-f(S)u-g(S)v\\
\dot{u}&=&(f(S)-D_0)u-\alpha(\cdot)u+\beta(\cdot)v\\
\dot{v}&=&(g(S)-D_1)v+\alpha(\cdot)u-\beta(\cdot)v
\end{array}
\right..
\end{equation}
In these equations, $S(t)$ denotes the concentration of the substrate at time $t$;
$u(t)$ and $v(t)$ denote, respectively, the concentration of the population of planktonic
microorganisms and flocks of bacteria at time $t$;
$f(S)$ and $g(S)$ represent, respectively, the per-capita growth rate of the populations;
$S_{in}$ and $D$ denote, respectively, the concentration of substrate in the feed
bottle and the dilution rate of the chemostat;
$D_0$ and $D_1$ represent, respectively, the removal rate of the microorganisms.

The dot in
attachment rate $\alpha(\cdot)$ and detachment rate $\beta(\cdot)$ means that these rates can depend
on the state variables, so that system (\ref{genmod}) recovers some  of
the models which were considered in the existing literature. For instance, the model of adaptive nutrient uptake, where
$u$ denotes the low growing cells and
$v$ denotes the fast growing cells considered in  \cite{tang} is obtained with attachment and detachment rates depending only on $S$
$$\alpha(\cdot)=\alpha(S),\qquad \beta(\cdot)=\beta(S).$$
The model of wall growth,
where $u$ denotes the density of planktonic bacteria, and $v$ denotes the density of wall-attached bacteria,
considered by Pilyugin and Waltman \cite{PW99}, is obtained with constant rates
$$\alpha(\cdot)=a,\qquad\beta(\cdot)=b.$$
The Freter model \cite{freter2} is given by
\begin{equation}\label{modFreter}
\left\{
\begin{array}{lcl}
\dot{S}&=&D(S_{in}-S)-f(S)u-g(S)v\\
\dot{u}&=&(f(S)-D_0)u-a(1-W)u+bv+g(S)(1-G(W))v\\
\dot{v}&=&(g(S)G(W)-D_1-b)v+a(1-W)u
\end{array}
\right.
\end{equation}
where $W=v/v_{{\rm max}}$ and $G(\cdot)$ is decreasing.
Notice that this model is a particular case of the model (\ref{genmod})
with
$$\alpha(\cdot)=a(1-W),\quad \beta(\cdot)=b+g(S)(1-G(W)) \ .$$
Actually, if $v_{{\rm max}}=\infty$ one obtains $W=0$ and if
 $G(0)=1$ then $\alpha(\cdot)=a$, $\beta(\cdot)=b$, and (\ref{modFreter}) is simply
the model of Pilyugin and Waltman  \cite{PW99}.

A model with flocks of two bacteria has been considered in \cite{HR08}
\begin{equation}\label{modHR}
\left\{
\begin{array}{lcl}
\dot{S}&=&D(S_{in}-S)-f(S)u-g(S)v\\
\dot{u}&=&(f(S)-D_0)u-au^2+bv\\
\dot{v}&=&(g(S)-D_1)v+au^2-bv
\end{array}
\right.
\end{equation}
that is a particular case of model (\ref{genmod}) obtained with
$$\alpha(\cdot)=au,\qquad\beta(\cdot)=b.$$
This model has been studied by Haegeman and Rapaport \cite{HR08} in the case of $g(\cdot)=0$ where the bacteria
in flocks are assumed to do not consume any substrate, and by Fekih-Salem and al \cite{Fekih} in the more general case of
$0\leqslant g(\cdot)\leqslant f(\cdot)$ where the bacteria in flocks consume less substrate than the isolated bacteria.
This model has been also extended to the case of
flocks with an arbitrary numbers of bacteria in \cite{HLH07}.\\

In the present paper we will not consider the size or the number of bacteria in flocks in our model.
We simply distinguish the biomass in flocks and the isolated biomass.
The biomass of isolated bacteria is denoted by $u$ and the biomass in flocks is denoted by $v$.
Hence isolated bacteria and isolated bacteria or flocks can stick together
to form new flocks,  with rate $a(u+v)u$, where $a$ is a constant, proportional to both the density of isolated bacteria,
that is $u$, and the total biomass density, that is $u+v$,
and flocks can split and liberate isolated bacteria,
with rate $bv$, where $b$ is a constant, proportional the their density $v$.
Hence, taking
$$\alpha(\cdot)=a(u+v),\qquad\beta(\cdot)=b$$
in model (\ref{genmod}), one obtains the following dynamical system
\begin{equation}\label{modSLRH}
\left\{
\begin{array}{l}
\dot{S}=D(S_{in}-S)-f(S)u-g(S)v\\
\dot{u}=(f(S)-D_0)u-a(u+v)u+bv\\
\dot{v}=(g(S)-D_1)v+a(u+v)u-bv
\end{array}
\right.
\end{equation}
that we aim to study in the present paper, as well as its extensions
to multi-species populations.

\section{A two time scales dynamics}
The general model for the flocculation is
\begin{equation*}
\left\{
\begin{array}{lcl}
\dot{S}&=&D(S_{in}-S)-f(S)u-g(S)v\\
\dot{u}&=&(f(S)-D_0)u
-\alpha(S,v,u)u+\beta(S,u,v)v\\
\dot{v}&=&(g(S)-D_1)v
+\alpha(S,v,u)u-\beta(S,u,v)v
\end{array}
\right..
\end{equation*}
If we assume that the dynamics of flocculation and deflocculation is much faster than the growth of the species, one can write the model the following way
\begin{equation}\label{sf1}
\left\{
\begin{array}{lcl}
\dot{S}&=&D(S_{in}-S)-f(S)u-g(S)v\\
\dot{u}&=&\ds (f(S)-D_0)u-\frac{\alpha(S,u,v)}{\varepsilon}u+\frac{\beta(S,u,v)}{\varepsilon}v\\
\dot{v}&=&\ds (g(S)-D_1)v+\frac{\alpha(S,u,v)}{\varepsilon}u-\frac{\beta(S,u,v)}{\varepsilon}v
\end{array}
\right.
\end{equation}
where $\varepsilon$ is expected to be a small non-negative number.
Notice that the dynamics of the total biomass $x=u+v$ is given by the equation
$$\dot{x}=(f(S)-D_0)u+(g(S)-D_1)v.$$
Thus, $u$ and $v$ are fast variables, while $S$ and $x$ are slow ones.
The fast dynamics is given by
$$u'=-{\alpha(S,u,v)}u+{\beta(S,u,v)}v \ , $$
and the slow manifold is defined by the positive solutions of the system
$$\alpha(S,u,v)u=\beta(S,u,v)v \quad \mbox{with} \quad u+v=x \ .$$
Hence one has
$$u=p(S,x)x,\qquad v=(1-p(S,x))x.$$
Assuming that this slow manifold is asymptotically stable for the fast
equation, the reduction
of the system to the slow system gives an approximation of the
solutions of (\ref{sf1}) for arbitrarily small $\varepsilon$, in
accordance with the theories of singular perturbations or variables aggregation \cite{O91,AP98,SAP11}.
The reduced model is obtained by replacing the fast variables $u$ and $v$
in the equations of $S$ and $x$:
\begin{equation}
\label{reduced-general}
\left\{
\begin{array}{l}
\dot{S}=D(S_{in}-S)-\mu(S,x)x\\
\dot{x}=\left(\mu(S,x)-d(S,x)\right)x
\end{array}
\right.
\end{equation}
{where}
$$
\mu(S,x)=p(S,x)f(S)+(1-p(S,x))g(S),
$$
$$
d(S,x)=p(S,x)D_0+(1-p(S,x))D_1.
$$
{Notice that $p(\cdot)$ depends on functions $\alpha(.)$ and $\beta(.)$}.
Consequently the growth function $\mu(\cdot)$ is a density dependent growth function model, as already studied in \cite{LH06,LMR04,LMR06}.
But, here the removal rates are replaced by a function $d(\cdot)$ that depends also on functions $\alpha(.)$ and $\beta(.)$. This last
property is new in the literature, to our knowledge.
\medskip

For the slow and fast case of (\ref{modHR}) one has, see \cite{HR08}
$$
\frac{\alpha(\cdot)}{\varepsilon}=\frac{a}{\varepsilon}u,\quad
\frac{\beta(\cdot)}{\varepsilon}=\frac{b}{\varepsilon}\quad
\mbox{and}\quad
p(x)=\frac{2}{1+\sqrt{1+4\frac{a}{b}x}}.$$
The slow and fast case of (\ref{modSLRH}) is given by
\begin{equation}\label{modSLRHsf}
\left\{
\begin{array}{l}
\dot{S}=D(S_{in}-S)-f(S)u-g(S)v\\
\dot{u}=(f(S)-D_0)u-\frac{a}{\varepsilon}(u+v)u+\frac{b}{\varepsilon}v\\
\dot{v}=(g(S)-D_1)v+\frac{a}{\varepsilon}(u+v)u-\frac{b}{\varepsilon}v
\end{array}
\right..
\end{equation}
The slow manifold is defined by the positive solutions of the system
$$a(u+v)u=bv \quad \mbox{with} \quad u+v=x \ .$$
Hence one has
$$u=p(x)x,\quad v=(1-p(x))x,\qquad\mbox{where}\quad
p(x)=\frac{b}{b+ax}.$$
Hence, the reduced model corresponding to (\ref{modSLRHsf}) is given by
\begin{equation}\label{reduced}
\left\{
\begin{array}{l}
\dot{S}=D(S_{in}-S)-\mu(S,x)x\\
\dot{x}=\left(\mu(S,x)-d(x)\right)x
\end{array}
\right.
\end{equation}
where
\begin{equation}\label{reduced1}
\mu(S,x)=\frac{bf(S)+axg(S)}{b+ax},\qquad
d(x)=\frac{bD_0+axD_1}{b+ax}.
\end{equation}
In the rest of the section, we use Tikhonov's theory \cite{LST,T,W}
(see also \cite{K02}),
and we justify that the solutions of (\ref{modSLRHsf}) are approximated by the solutions of the reduced model (\ref{reduced},\ref{reduced1}).
\begin{theorem}
Let $(S(t,\varepsilon),u(t,\varepsilon),v(t,\varepsilon))$ be the solution of (\ref{modSLRHsf}) with initial condition $(S_0,u_0,v_0)$ satisfying
$S_0\geqslant 0$, $u_0>0$, and $v_0\geqslant0$.
Let $(\overline{S}(t),\overline{x}(t))$ be the solution of the reduced problem (\ref{reduced})
with initial conditions
$$\overline{S}(0)=S_0,\quad \overline{x}(t)=u_0+v_0.$$
Then as $\varepsilon\to 0$,
\begin{equation}\label{approxSx}
S(t,\varepsilon)=\overline{S}(t)+o(1),\qquad x(t,\varepsilon)=\overline{x}(t)+o(1)
\end{equation}
uniformly for $t\in[0,T]$, and
\begin{equation}\label{approxuv}
u(t,\varepsilon)=\frac{b}{b+a\overline{x}(t)}\overline{x}(t)+o(1),\qquad
v(t,\varepsilon)= \frac{a\overline{x}(t)}{b+a\overline{x}(t)}\overline{x}(t)+o(1)
\end{equation}
uniformly for $t\in[t_0,T]$, where $0<t_0<T$ are arbitrary but fixed and independent of $\varepsilon$.
If the solution $(\overline{S}(t),\overline{x}(t))$  of the reduced  problem converges to an asymptotically stable equilibrium, then we can put $T=+\infty$ in the the approximations (\ref{approxSx}) and (\ref{approxuv}) given.
\end{theorem}

\noindent
{\bf Proof}.
In the variables $S$, $x=u+v$ and $u$, system (\ref{modSLRHsf}) is written
\begin{equation}\label{eqent2}
\left\{
\begin{array}{l}
\dot{S}=D(S_{in}-S)-f(S)u-g(S)(x-u)\\
\dot{x}=f(S)u+g(S)(x-u)-D_0 u-D_1(x-u)\\
\dot{u}=\displaystyle(f(S)-D_0)u-\frac{a}{\varepsilon}xu+\frac{b}{\varepsilon}(x-u)
\end{array}
\right..
\end{equation}
This is a slow and fast system with two slow variables $S$ and $x$ and one fast variable, $u$.
The fast equation is
\begin{equation}\label{fe}
{u}'=-{ax}u+{b}(x-u)=bx-(ax+b)u.
\end{equation}
The slow manifold (or quasi-steady-state) is given by
$$u=xp(x),\quad \mbox{ where }\quad p(x)=\frac{b}{b+ax}.$$
Since $ax+b>0$ this slow manifold is globally asymptotically
stable. Thus, the Tikhonov's Theorem \cite{LST,T,W}
applies and asserts that
after a fast transition toward the slow manifold, the solutions are approximated by a solution of the reduced equation,
which is obtained by replacing the fast variable $u$ in (\ref{eqent2}) by the quasi steady state $u=xp(x)$. This reduced system is
\begin{equation*}
\left\{
\begin{array}{l}
\dot{S}=D(S_{in}-S)-\left(f(S)p(x)+g(S)(1-p(x))\right)x\\
\dot{x}=\left(f(S)p(x)+g(S)(1-p(x))-D_0p(x)-D_1(1-p(x))\right)x
\end{array}
\right..
\end{equation*}
This is system (\ref{reduced},\ref{reduced1}). The approximations
(\ref{approxSx}) and (\ref{approxuv}) follow from the Tikhonov's  Theorem. Recall that when the reduced problem has an asymptotically stable equilibrium,
then these approximations hold for all $t>0$ and not only on a compact interval $[0,T)$. Recall also that there is a boundary layer for the fast
variables $u$ and $v$, that is the approximations (\ref{approxuv}) hold only for $t\geqslant t_0$ where $t_0>0$ can be arbitrarily small but fixed.
Therefore, very quickly, the density $v(t)$ of flocks tends to
$$\frac{a({x}(0))^2}{b+{a}{x}(0)},$$
and the density $u(t)$ of planktonic bacteria tends toward
$$\frac{b{x}(0)}{b+{a}{x}(0)}.$$
These values depend only on the initial total biomass $x(0)=u(0)+v(0)$, and on the ratio
$a/b$ between the rate of flocculation and deflocculation. After this boundary layer, one has a slow variation
of the densities of flocks and isolated bacteria according to the density dependent reduced model (\ref{reduced},\ref{reduced1}).

\section{Study of the reduced model for one species}
\label{bistability}
We consider in this section the mathematical analysis of the reduced model (\ref{reduced}). We do not assume that the functions $\mu(S,x)$ and $d(x)$ are of the particular form (\ref{reduced1}).
We assume that
\begin{description}
\item[H0:] $\mu(0,x)=0$ and $\mu(S,x)>0$ for all $S>0$ and all
$x \geqslant 0$.
\item[H1:] $\displaystyle\frac{\partial{\mu}}{\partial S}>0$ and
$\displaystyle\frac{\partial{\mu}}{\partial x}<0$ for all $S>0$ and all $x\geqslant 0$.
\item[H2:] $d(0)=D_0$, $d(+\infty)=D_1<D_0\leqslant D$, $d(x)>0$, $d'(x)<0$ and $[xd(x)]'>0$ for all $x\geqslant 0$.
\end{description}
Let us denote by
$$f(S)=\mu(S,0),\mbox{ and } g(S)=\mu(S,+\infty).$$
The functions $f(\cdot)$ and $g(\cdot)$ are increasing and positive for all $S>0$.
If equations $f(S)=D_0$ and $g(S)=D_1$ have solutions, one lets
$$\lambda_0=f^{-1}(D_0),\mbox{ and }\lambda_1=g^{-1}(D_1).$$
Otherwise one lets $\lambda_k=+\infty$, $k=0,1$.
We add the following assumption
\begin{description}
\item[H3:] If $\lambda_0<\lambda_1$, then for all $S\in[\lambda_0,\lambda_1)$ and $x\geqslant 0$ one has
$
\displaystyle d'(x)>\frac{\partial{\mu}}{\partial x}(S,x).
$
\item[H4:] If $\lambda_1<\lambda_0$, then for all $S\in(\lambda_1,\lambda_0]$ and $x\geqslant 0$ one has
$
\displaystyle d'(x)<\frac{\partial{\mu}}{\partial x}(S,x).
$
\end{description}
\begin{lemma}
Assumptions {\bf H0}-{\bf H4} are satisfied in the case
$$
{\mu}(S,x)=f(S)p(x)+g(S)(1-p(x)),
\quad
d(x)=D_{0}p(x)+D_1(1-p(x))
$$
where $f(\cdot)>g(\cdot)$ are increasing functions and
$p(x)$ is a decreasing functions such that
$p(0)=1$, $p(+\infty)=0$ with $[xp(x)]'>0$.
\end{lemma}

\begin{figure}[ht]
\setlength{\unitlength}{1.0cm}
\begin{center}
\begin{picture}(9,4)(0,0)
\put(-0.5,-0.2){\rotatebox{0}{\includegraphics[scale=0.2]{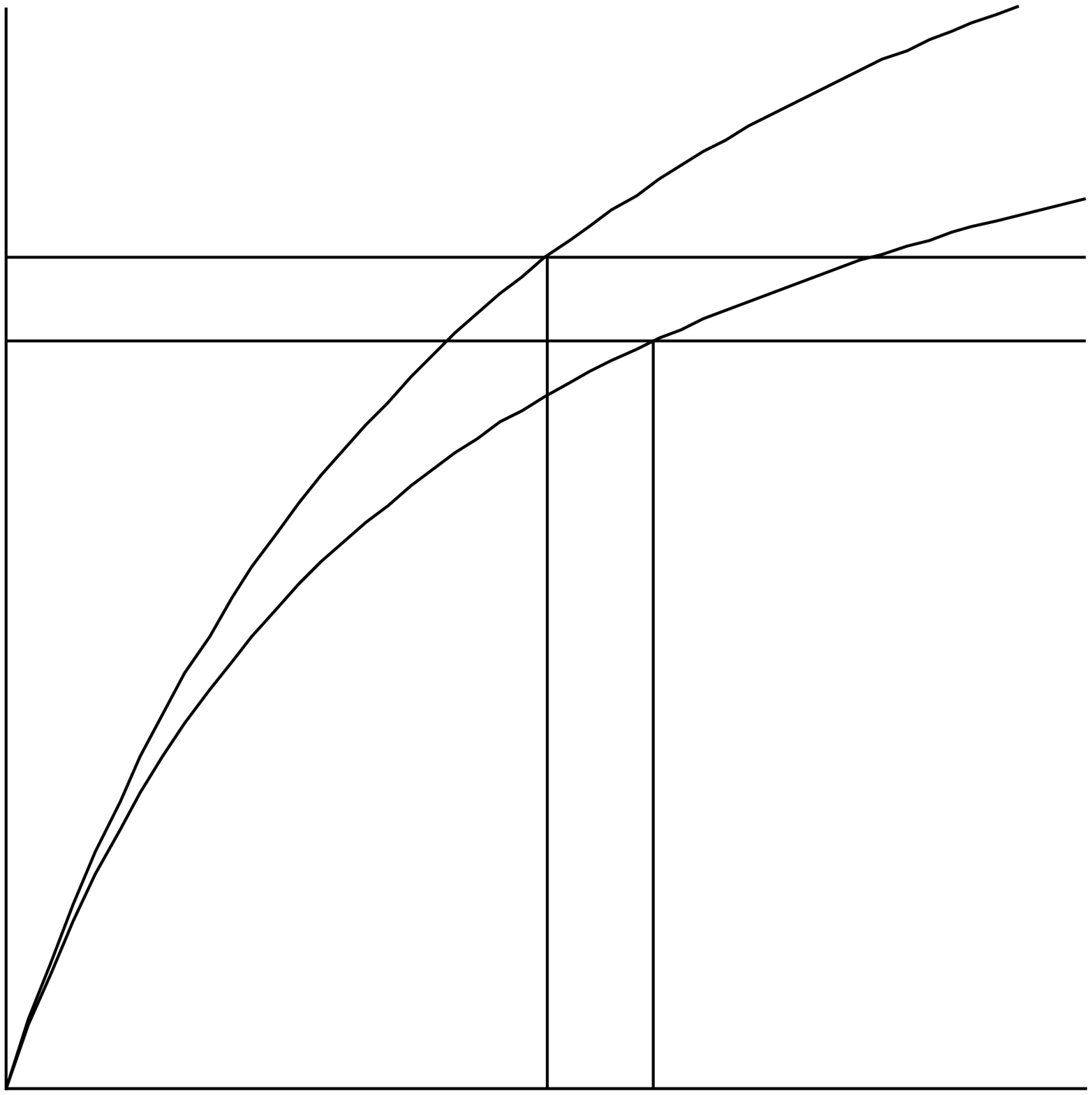}}}
\put(5.5,-0.2){\rotatebox{0}{\includegraphics[scale=0.2]{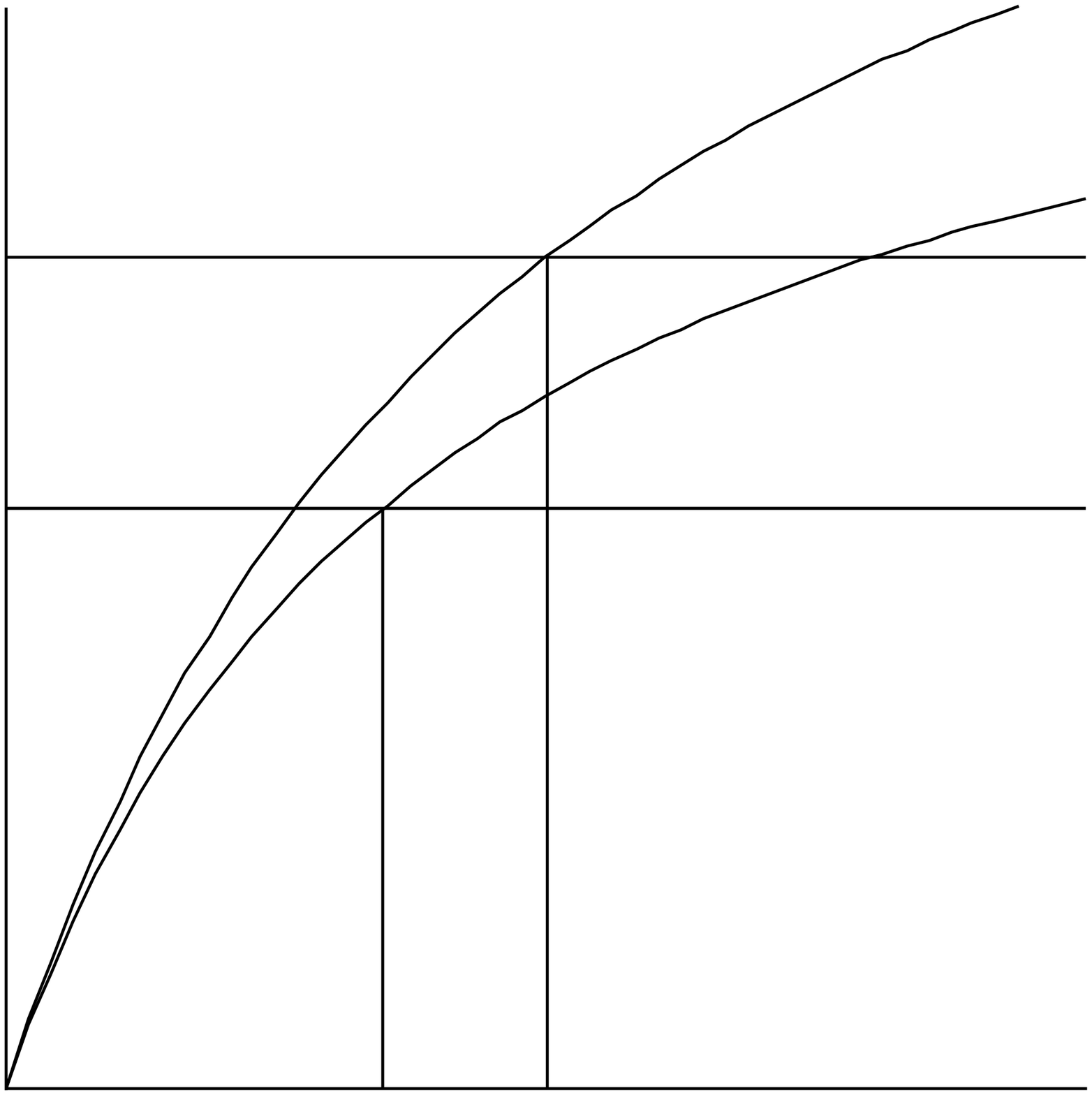}}}
\put(3.3,3.6){{\small $f(S)$}}
\put(3.5,3){{\small $g(S)$}}
\put(9.3,3.6){{\small $f(S)$}}
\put(9.5,3){{\small $g(S)$}}
\put(-0.9,2.7){{\small $D_0$}}
\put(-0.9,2.4){{\small $D_1$}}
\put(5.1,2.7){{\small $D_0$}}
\put(5.1,1.8){{\small $D_1$}}
\put(1.3,-0.5){{\small $\lambda_0$}}
\put(1.8,-0.5){{\small $\lambda_1$}}
\put(7.4,-0.5){{\small $\lambda_0$}}
\put(6.8,-0.5){{\small $\lambda_1$}}
\put(3.3,-0.4){{\small $S$}}
\put(9.3,-0.4){{\small $S$}}
\end{picture}
\end{center}
\caption{On the left, the case $\lambda_0<\lambda_1$.
On the right, the case $\lambda_1<\lambda_0$.}\label{lesmu}
\end{figure}

\noindent
{\bf Proof}.
Since $p'(x)<0$ and $f(S)>g(S)$ for all $S>0$, one has
$$
\begin{array}{l}
\displaystyle \frac{\partial{\mu}}{\partial S}=f'(S)p(x)+g'(S)(1-p(x))>0,\\[3mm]
\displaystyle \frac{\partial{\mu}}{\partial x}=(f(S)-g(S))p'(x)<0.
\end{array}$$
Thus {\bf H1} is satisfied.
On the other hand
$$d'(x)=(D_0-D_1)p'(x)<0,\qquad [xd(x)]'=D_1+(D_0-D_1)[xp(x)]'>0$$
since $p'(x)<0$, $[xp(x)]'>0$ and $D_{0}>D_{1}$.
Thus {\bf H2} is satisfied. Moreover
$$\frac{\partial{\mu}}{\partial x}(S,x)-d'(x)=[f(S)-D_0+D_1-g(S)]p'(x)<0$$
for all $S\in[\lambda_0,\lambda_1)$ and $x\geqslant 0$, since $p'(x)<0$ and
$f(S)\geqslant D_0$,  $D_1>g(S)$ for $\lambda_0\leqslant S<\lambda_1$, see Figure \ref{lesmu}.
Thus Assumption  {\bf H3} is also satisfied.
Similarly one has
$$\frac{\partial{\mu}}{\partial x}(S,x)-d'(x)=[f(S)-D_0+D_1-g(S)]p'(x)>0$$
for all $S\in(\lambda_1,\lambda_0]$ and $x\geqslant 0$, since $p'(x)<0$ and
$f(S)\leqslant D_0$,  $D_1<g(S)$ for $\lambda_1<S\leqslant\lambda_0$, see Figure \ref{lesmu}.
Thus Assumption {\bf H4} is satisfied.\\

\medskip
\noindent
Notice that, if
$p(x)=\frac{b}{b+ax},$
then the properties of $p(\cdot)$ stated in the lemma are satisfied. Indeed
$$
p'(x)=\frac{-ab}{(b+ax)^2},\qquad [xp(x)]'=\frac{b^2}{(b+ax)^2}.
$$

\subsection{Existence of equilibria}
The equilibria of the system are solutions of the set of equations
$$
\left\{
\begin{array}{l}
D(S_{in}-S)-\mu(S,x)x=0\\
(\mu(S,x)-d(x))x=0
\end{array}
\right.
$$
The second equation is equivalent to $x=0$ or $\mu(S,x)=d(x)$.
If $x=0$ then from the first equation one has $S=S_{in}$.
This is the washout equilibrium
$$E_0=(S_{in},0).$$
If $\mu(S,x)=d(x)$, the first equation gives $D(S_{in}-S)=xd(x)$.
Hence
$$
S=\gamma(x):=S_{in}-\frac{xd(x)}{D}.
$$
Since
$$\gamma(0)=S_{in},\mbox{ and }\gamma'(x)=-\frac{[xd(x)]'}{D}<0$$
the function $\gamma(\cdot)$ is decreasing.
Thus one have to solve the equation
\begin{equation}\label{equation}
\mu(S,x)=d(x).
\end{equation}
Since $\frac{\partial{\mu}}{\partial S}>0$,
by the implicit function theorem, thus equation defines a function
$$S=\phi(x),\mbox{ such that }\lambda_0=\phi(0),$$
and
$$\phi'(x)=\frac{\displaystyle d'(x)-\frac{\partial{\mu}}{\partial x}(\phi(x),x)}{\displaystyle \frac{\partial{\mu}}{\partial S}(\phi(x),x)}.$$
The sign of $\phi'(\cdot)$ is given by assumptions {\bf H1}, {\bf H3} and {\bf H4}. The cases $\lambda_0<\lambda_1$ and $\lambda_0>\lambda_1$ have to be distinguished.

\begin{figure}[ht]
\setlength{\unitlength}{1.0cm}
\begin{center}
\begin{picture}(9,4)(0,0)
\put(0,-0.2){\rotatebox{0}{\includegraphics[scale=0.2]{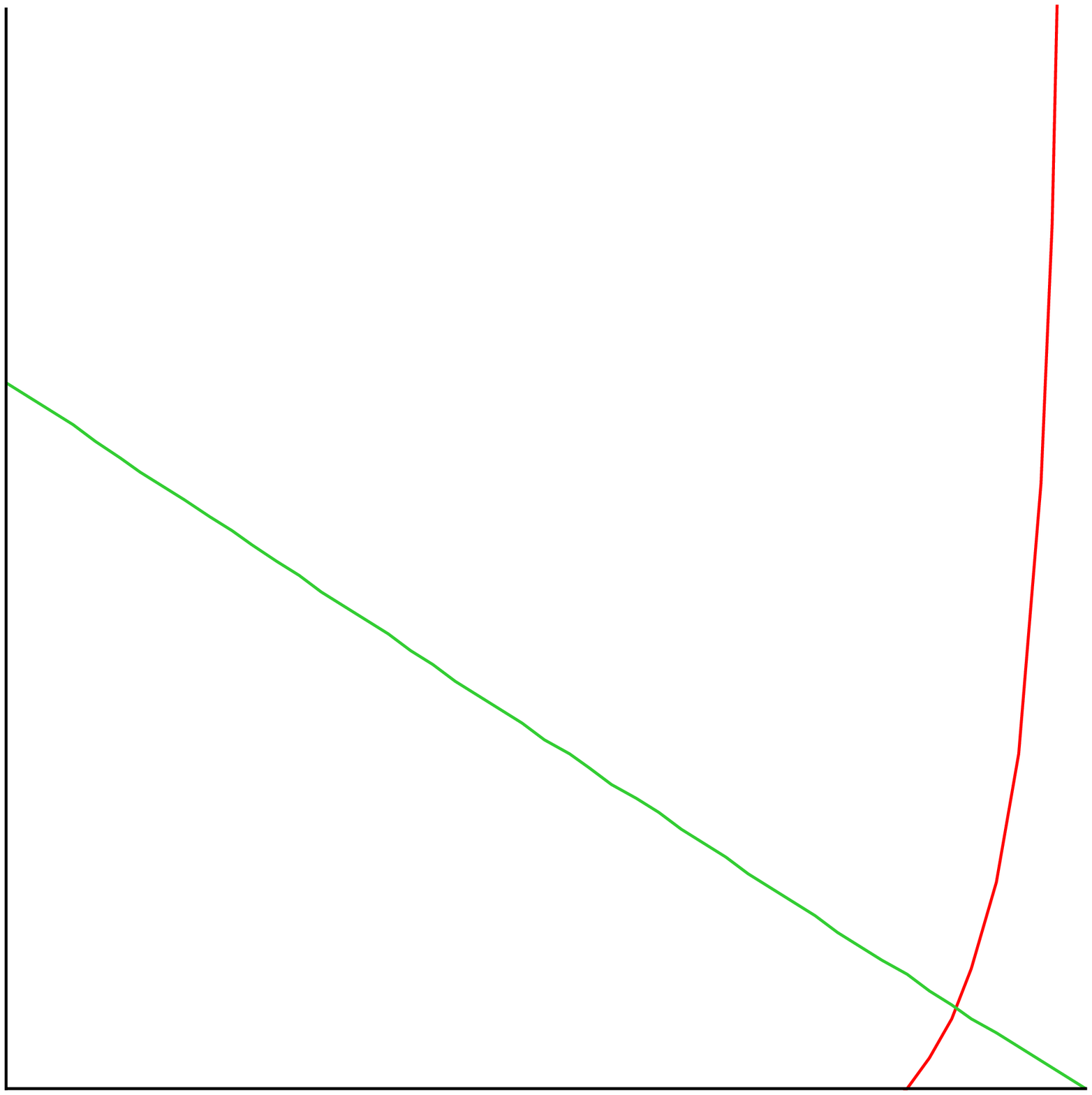}}}
\put(5,-0.2){\rotatebox{0}{\includegraphics[scale=0.2]{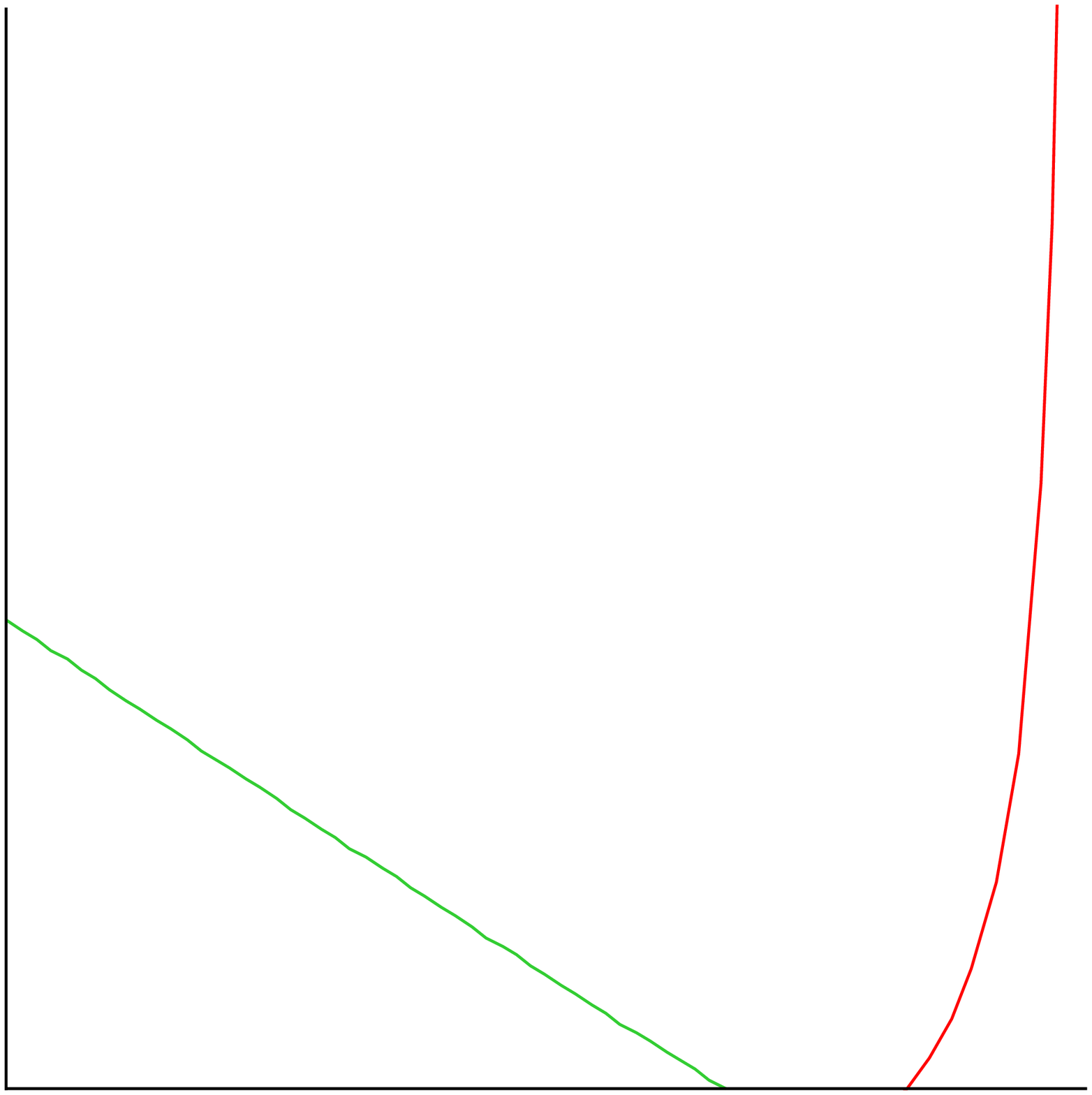}}}
\put(3.6,-0.5){{\small $S_{in}$}}
\put(2.4,3.3){{\small $S=\phi(x)$}}
\put(1,1.8){{\small $S=\gamma(x)$}}
\put(3,-0.5){{\small $\lambda_0$}}
\put(7.4,3.3){{\small $S=\phi(x)$}}
\put(6,1){{\small $S=\gamma(x)$}}
\put(7.4,-0.5){{\small $S_{in}$}}
\put(8.1,-0.5){{\small $\lambda_0$}}
\put(4,-0.2){{\small $S$}}
\put(0.2,3.5){{\small $x$}}
\put(9,-0.2){{\small $S$}}
\put(5.2,3.5){{\small $x$}}
\end{picture}
\end{center}
\caption{Null-clines $S=\phi(x)$ and $S=\gamma(x)$ in the case $\lambda_0<\lambda_1$. On the left, the case $\lambda_0<S_{in}$ with a unique intersection point. On the right, the case $\lambda_0>S_{in}$ with no intersection point.} \label{figcas1}
\end{figure}

When $\lambda_0<\lambda_1$
the function $S=\phi(x)$ is defined for all $x\geqslant 0$ and satisfies
$$\lambda_0=\phi(0),\qquad \lambda_1=\phi(+\infty),\qquad \phi'(x)>0.$$
The equilibria are the intersection points of the graphs of functions
$$S=\phi(x)\mbox{ and }S=\gamma(x).$$
Since the first function is increasing and the second one is decreasing, there is a unique solution if
$\lambda_0<S_{in}$,
and no solution if $\lambda_0>S_{in}$, see Figure \ref{figcas1}.
\begin{proposition}
If $\lambda_0<\min(\lambda_1,S_{in})$, there exists a unique positive equilibrium.
If $S_{in}<\lambda_0<\lambda_1$, there is no positive equilibrium.
\end{proposition}

\begin{figure}[ht]
\setlength{\unitlength}{1.0cm}
\begin{center}
\begin{picture}(9,4)(0,0)
\put(0,-0.2){\rotatebox{0}{\includegraphics[scale=0.2]{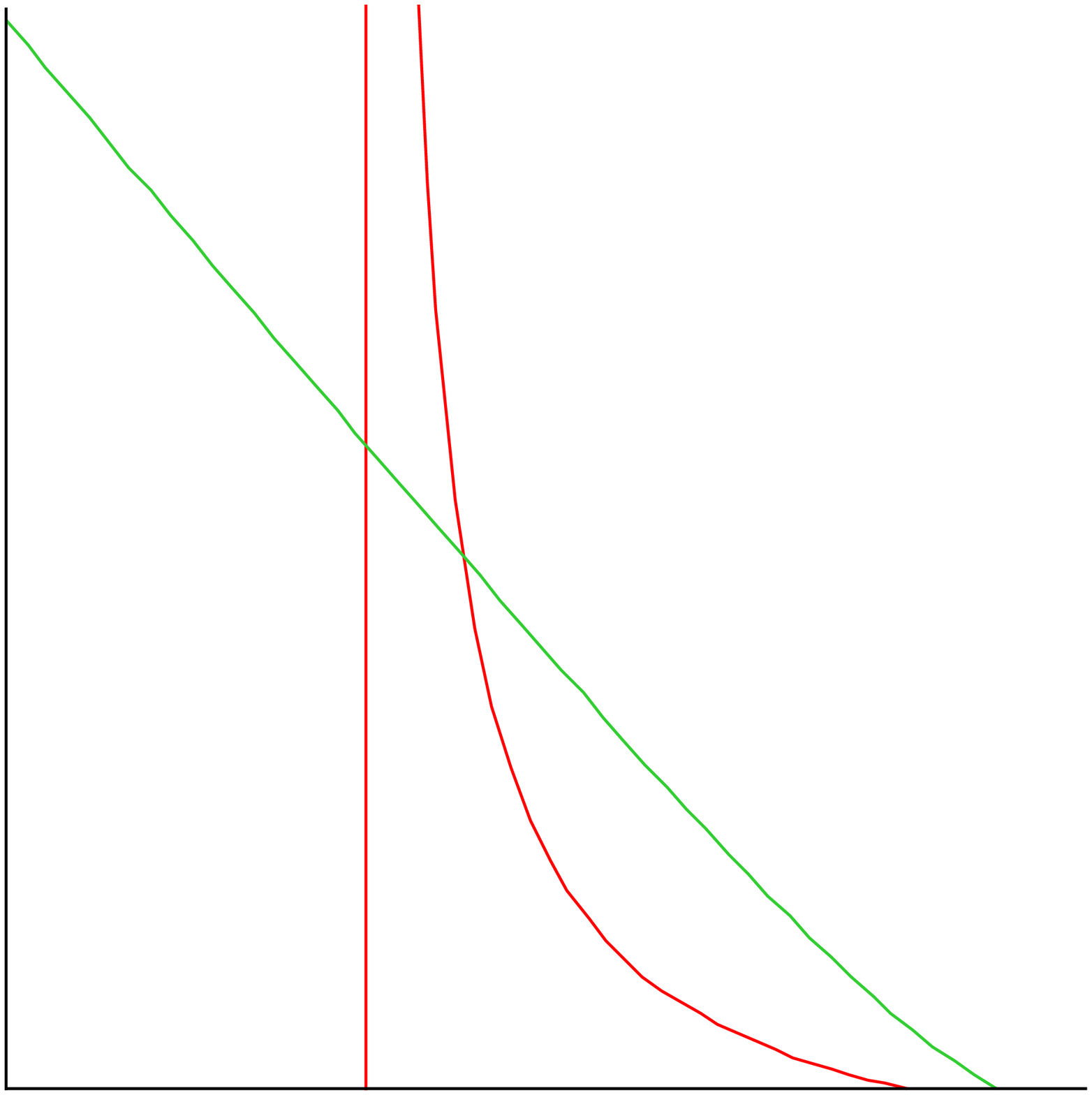}}}
\put(5,-0.2){\rotatebox{0}{\includegraphics[scale=0.2]{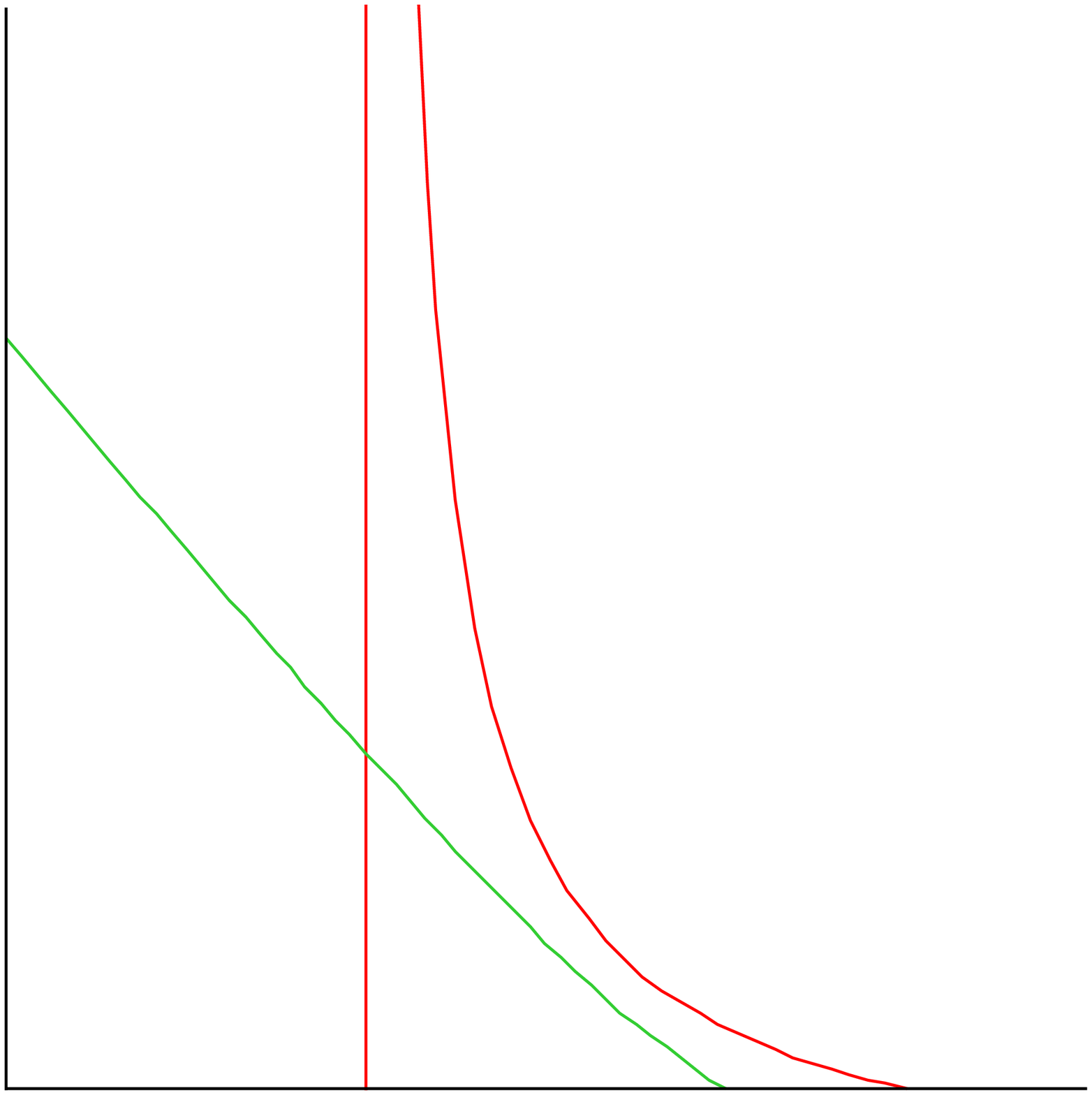}}}
\put(3.5,-0.5){{\small $S_{in}$}}
\put(1.6,3.3){{\small $S=\phi(x)$}}
\put(2.4,1){{\small $S=\gamma(x)$}}
\put(3,-0.5){{\small $\lambda_0$}}
\put(1.2,-0.5){{\small $\lambda_1$}}
\put(6.6,3.3){{\small $S=\phi(x)$}}
\put(5.1,2){{\small $S=\gamma(x)$}}
\put(7.4,-0.5){{\small $S_{in}$}}
\put(8.1,-0.5){{\small $\lambda_0$}}
\put(6.2,-0.5){{\small $\lambda_1$}}
\end{picture}
\end{center}
\caption{Null-clines $S=\phi(x)$ and $S=\gamma(x)$ in the case $\lambda_0>\lambda_1$. On the left, the case $\lambda_0<S_{in}$ with
at least one intersection point. On the right, the case $\lambda_0>S_{in}$ with no intersection point.
In the second case it is possible to have two intersection points, see Figure \ref{fig2}.} \label{figcas2}
\end{figure}

\begin{figure}[htbp]
\setlength{\unitlength}{1.0cm}
\begin{center}
\begin{picture}(11,5)(0,0.2)
\put(0,0){\rotatebox{0}{\includegraphics[scale=0.25]{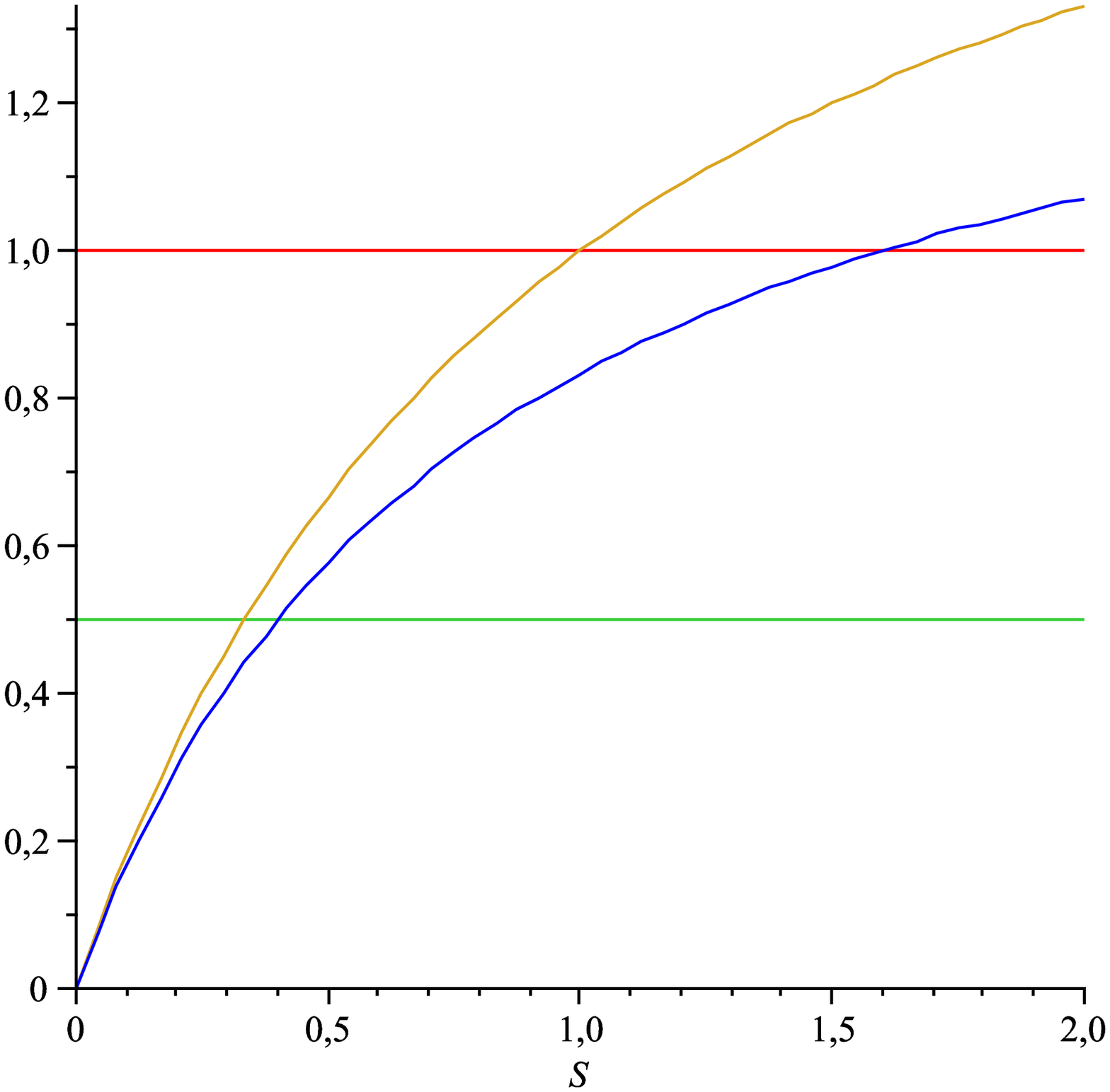}}}
\put(6,0){\rotatebox{0}{\includegraphics[scale=0.25]{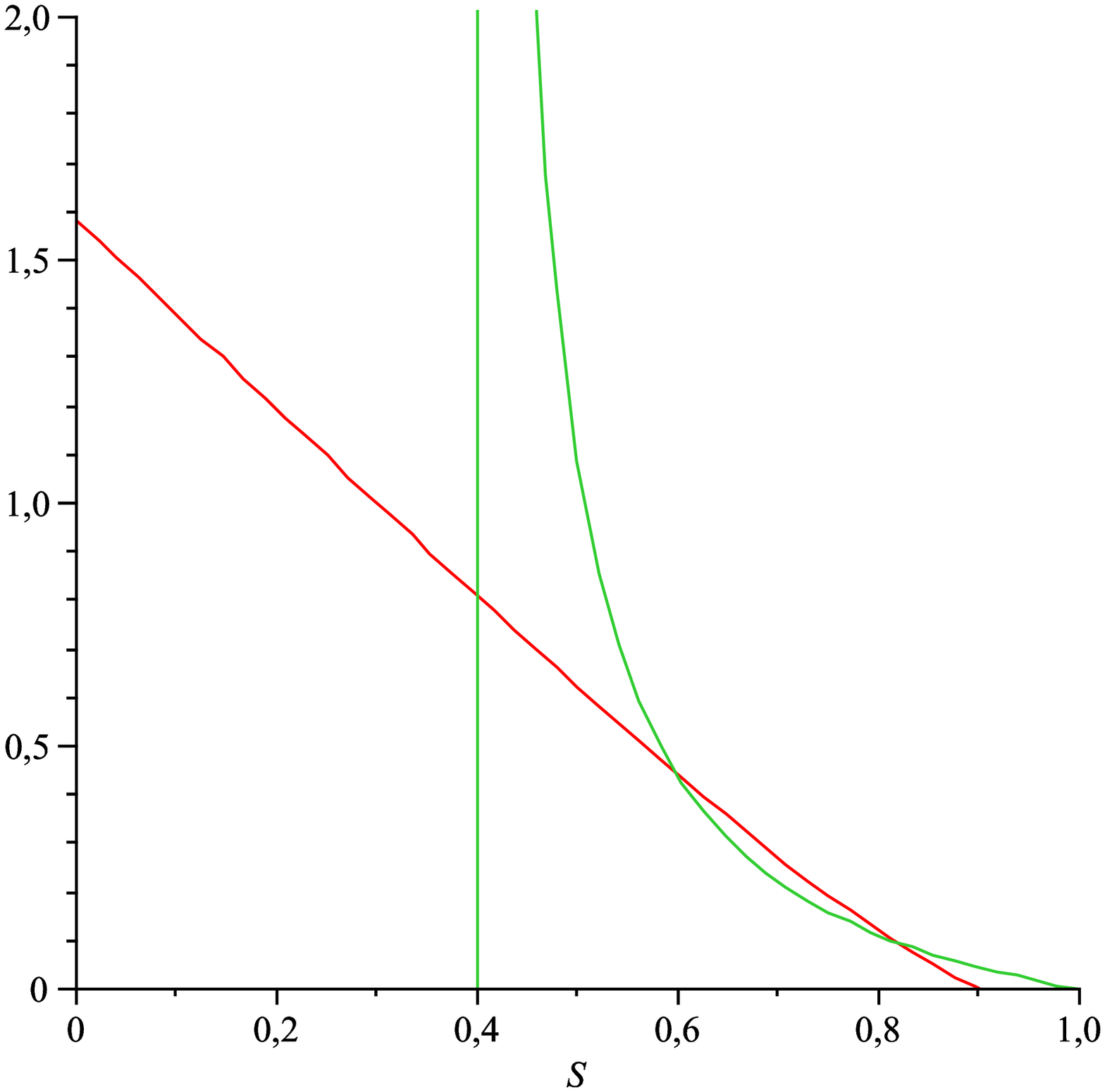}}}
\put(-0.4,3.6){{\small $D_0$}}
\put(-0.4,2){{\small $D_1$}}
\put(3.6,4.7){{\small $f(S)$}}
\put(4,4){{\small $g(S)$}}
\put(8.5,4.2){{\small $S=\phi(x)$}}
\put(6.7,3.6){{\small $S=\gamma(x)$}}
\put(8,0){{\small $\lambda_1$}}
\put(10,0){{\small $S_{in}$}}
\put(10.7,0){{\small $\lambda_0$}}
\end{picture}
\end{center}
\caption{On the left, the growth functions $f(\cdot)$ and $g(\cdot)$ given  by (\ref{parva1}).
On the right, the corresponding null-clines $S=\phi(x)$ and $S=\gamma(x)$ for the parameter values(\ref{parva2}),
showing the existence of two positive equilibria in the case when $\lambda_0>S_{in}$.} \label{fig2}
\end{figure}

When $\lambda_1<\lambda_0$
the function $S=\phi(x)$ is defined for all $x\geqslant 0$ and satisfies
$$\lambda_0=\phi(0),\qquad \lambda_1=\phi(+\infty),\qquad \phi'(x)<0.$$
Both functions
$S=\phi(x)$ and $S=\gamma(x)$ are decreasing and
$$\phi(0)=\lambda_0,\qquad S_{in}=\gamma(0).$$
Thus, if $\lambda_0<S_{in}$ the graphs intersect at at least one
non-negative point (see Figure \ref{figcas2}).
If $\lambda_0>S_{in}$, the graphs of functions $S=\phi(x)$
and $S=\gamma(x)$ can interest or not
(see Figures \ref{figcas2} and \ref{fig2}).

\begin{proposition}
If $\lambda_1<\lambda_0<S_{in}$, then there exists at least one positive equilibrium.
Generically one has an odd number of positive equilibria.
If $\lambda_1<S_{in}<\lambda_0$, then the system has generically no positive equilibrium or an even number of positive equilibria.
\end{proposition}

\subsection{Stability of equilibria}
The Jacobian matrix of (\ref{reduced}) is given by
$$
J=
\left[
\begin{array}{cc}
\displaystyle-D-x\frac{\partial{\mu}}{\partial S}(S,x)
&
\displaystyle-x\frac{\partial{\mu}}{\partial x}(S,x)-\mu(S,x)
\\[3mm]
\displaystyle x\frac{\partial{\mu}}{\partial S}(S,x)
&
\displaystyle\mu(S,x)-d(x)+x\frac{\partial{\mu}}{\partial x}(S,x)-xd'(x)
\end{array}
\right]
$$
At washout $E_0=(S_{in},0)$ this matrix is
$$
J_0=
\left[
\begin{array}{cc}
\displaystyle-D
&
\displaystyle-f(S_{in})
\\[3mm]
0
&
\displaystyle f(S_{in})-D_{0}
\end{array}
\right]
$$
The eigenvalues are $-D$ and $f(S_{in})-D_{0}$. Hence one has the following result

\begin{proposition}
If $\lambda_0<S_{in}$, then $E_0$ is unstable (saddle point). If
$\lambda_0>S_{in}$, then $E_0$ is locally exponentially stable (stable node).
\end{proposition}

At a positive equilibrium $E_1=(S,x)$, one has necessarily $\mu(S,x)=d(x)$, and the Jacobian matrix is
$$
J_1=
\left[
\begin{array}{cc}
\displaystyle-D-x\frac{\partial{\mu}}{\partial S}(S,x)
&
\displaystyle-x\frac{\partial{\mu}}{\partial x}(S,x)-\mu(S,x)
\\[3mm]
\displaystyle x\frac{\partial{\mu}}{\partial S}(S,x)
&
\displaystyle x\frac{\partial{\mu}}{\partial x}(S,x)-xd'(x)
\end{array}
\right]
$$
The trace of $J_1$ is
$${\rm tr}(J_1)=-D-x\frac{\partial{\mu}}{\partial S}(S,x)
+x\frac{\partial{\mu}}{\partial x}(S,x)-xd'(x).$$
Since
$[xd(x)]'=d(x)+xd'(x)$, one has
$$
{\rm tr}(J_1)=-D+d(x)-x\frac{\partial{\mu}}{\partial S}(S,x)
+x\frac{\partial{\mu}}{\partial x}(S,x)-[xd(x)]'<0
$$
since  $d(x)\leqslant D$, $\frac{\partial{\mu}}{\partial S}>0$, $\frac{\partial{\mu}}{\partial x}<0$
and $[xd(x)]'>0$.
The determinant of $J_1$ is
$$
{\rm det}(J_1)=Dx\left(
d'(x)-\frac{\partial{\mu}}{\partial x}\right)+
x\frac{\partial{\mu}}{\partial S}[xd(x)]'.
$$
One can write this determinant as
$$
{\rm det}(J_1)=
Dx\frac{\partial{\mu}}{\partial S}\left(
\frac{\displaystyle d'(x)-\frac{\partial{\mu}}{\partial x}}{\displaystyle \frac{\partial{\mu}}{\partial S}}
+\frac{[xd(x)]'}{D}\right)
=Dx\frac{\partial{\mu}}{\partial S}\left(
\phi'(x)-\gamma'(x)\right).
$$
Thus, If $\phi'(x)>\gamma'(x)$, then the determinant is positive and hence the eigenvalues are of negative real parts: the equilibrium $E_1$ is locally asymptotically stable.
On the other hand, if $\phi'(x)<\gamma'(x)$, then the determinant is negative, and the eigenvalues are of opposite sign:
the positive equilibrium is a saddle point. Thus we have shown the following result

\begin{proposition}
The following cases occur.
\begin{enumerate}
\item
If $\lambda_0<\min(\lambda_1,S_{in})$, the unique positive equilibrium is locally asymptotically stable (and fulfills $\phi'(x)>0>\gamma'(x)$).
\item If $\lambda_1<\lambda_0$ there are possibilities of multiple positive equilibria, that are saddle points when  $\phi'(x)<\gamma'(x)$ or stable nodes when  $\phi'(x)>\gamma'(x)$:
\begin{itemize}
\item[-] when $\lambda_0<S_{in}$, there exists at least one positive equilibrium and one has an odd number of equilibria which are alternatively stable and unstable,
\item[-]  when $\lambda_0>S_{in}$, there is no or an even  number of equilibria which are alternatively stable and unstable.
\end{itemize}
\end{enumerate}
\end{proposition}

\begin{figure}[htbp]
\setlength{\unitlength}{1.0cm}
\begin{center}
\begin{picture}(11,5)(0,0.2)
\put(3,0){\rotatebox{0}{\includegraphics[scale=0.25]{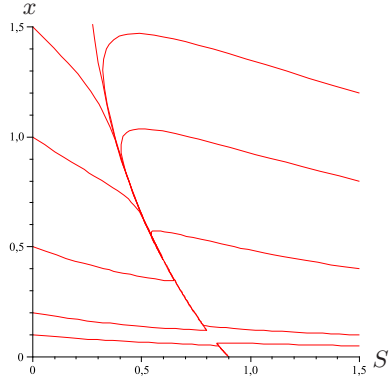}}}
\put(3.25,4.9){{\small $x$}}
\put(7.9,0.2){{\small $S$}}
\end{picture}
\end{center}
\caption{For the Monod functions (\ref{parva1}) and the parameters values (\ref{parva2}), the system  exhibits bi-stability.} \label{fig3}
\end{figure}

One concludes that in the case $\lambda_1<\lambda_0$, the system can exhibits bi-stability of the washout equilibrium and a positive equilibrium.
If the initial density of flocks is small enough, the solutions will converge to the washout equilibrium, otherwise
it will converge to the positive equilibrium, see Figure \ref{fig3}. In this figure we notice that the washout equilibrium is stable together with the positive equilibrium corresponding the lowest value of $S$. The domains of attraction of the stable equilibria are separated by the stable separatrix of the positive saddle node. The simulations shown in Figures  \ref{fig2} and \ref{fig3} where obtained for (\ref{reduced},\ref{reduced1}), with
the following Monod functions
\begin{equation}\label{parva1}
f(S)=\frac{2\,S}{1+S},\qquad g(S)=\frac{1.5\,S}{0.8+S}
\end{equation}
and the following values of the parameters
\begin{equation}\label{parva2}
D_0=D=1,\quad D_1=0.5,\quad a=4,\quad b=1,\quad S_{in}=0.9.
\end{equation}

\section{Flocculation with several species}
\label{multispecies}

We assume that $n$ species are competing on a same limiting resource,
and that each species is present in two forms:
isolated bacteria, of density $u_i$, and bacteria in flocks, of density $v_i$, for $i=1\cdots n$.
We assume that isolated bacteria can stick with isolated bacteria with
flocks to form new flocks, with rate $\alpha_i(\cdot)u_i$.
We assume also that flocks can split and liberate isolated bacteria
with rate $\beta_i(\cdot)v_i$.
$$u_i\xrightarrow{\alpha_i(\cdot)u_i}v_i,\qquad u_i\xleftarrow{\beta_i(\cdot)v_i}v_i \ .$$
Then the equations are
\begin{equation}\label{Nspecies}
\left\{
\begin{array}{l}
\dot{S}=\displaystyle D(S_{in}-S)-\sum_{i=1}^n(f_{i}(S)u_i+g_{i}(S)v_i)\\[3mm]
\dot{u}_i=(f_{i}(S)-D_{0i})u_i-\alpha_i(\cdot)u_i+\beta_i(\cdot)v_i, \quad 1\leqslant i\leqslant n\\
\dot{v}_i=(g_{i}(S)-D_{1i})v_i+\alpha_i(\cdot)u_i-\beta_i(\cdot)v_i
\end{array}
\right.
\end{equation}
The dynamics of the total biomass densities $x_i=u_i+v_i$
of the species $i$ is
$$\dot{x_i}=f_{i}(S)u_i+g_{i}(S)v_i-D_{0i}u_i-D_{1i}v_i \ .$$
We consider here the case where
$$
\alpha_i(\cdot)=\sum_{j=1}^nA_{ij}x_{j},\qquad \beta_i(\cdot)=B_i
$$
where $A_{ij}$ and $B_i$ are non-negative constants.

\begin{remark}
By letting the functions $\alpha_{i}(\cdot)$ depending on $x_{j}$ with $j\neq i$, we implicitly consider that flocks or aggregates can incorporate individuals of different species.\\
The removal rate $D_{1i}$ has to be interpreted as the sum
of the removal rate of the aggregates and the mortality rate of each
species $i$ in its attached form.
\end{remark}

Moreover, we assume that the dynamics of
flocculation and deflocculation are fast compared with the dynamics of the growth of bacteria, that is
$$
A_{ij}=\frac{a_{ij}}{\varepsilon},\quad B_i=\frac{b_i}{\varepsilon}
\ .
$$
In the variables $S$, $x_i$, $u_i$, the system is written
\begin{equation}\label{eqent3}
\left\{
\begin{array}{lcl}
\dot{S}&=&\ds D(S_{in}-S)-\sum_{i=1}^{n}[f_{i}(S)u_i+g_{i}(S)(x_{i}-u_{i})]\\[3mm]
\dot{x_i}&=&\ds f_{i}(S)u_i+g_{i}(S)(x_i-u_i)-D_{0i}u_i-D_{1i}(x_i-u_i), \quad i=1\cdots n\\
\dot{u}_i&=&\ds (f_{i}(S)-D_{0i})u_i-\frac{1}{\varepsilon}\sum_{j=1}^{n}a_{ij}x_ju_{i}+\frac{b_i}{\varepsilon}(x_i-u_i)
\end{array}
\right.
\end{equation}
This is a slow/fast system with the variables $S$, $x_i$ are
slow and the variables $u_i$ fast.
The fast equations are
\begin{equation}\label{fe3}
{u}_i'=-\sum_{j=1}^{n}a_{ij}x_ju_i+{b_i}(x_i-u_i),
\quad i=1\cdots n
\end{equation}
where $x_i$ are considered as parameters.
The slow manifold (or quasi steady-state) is given by
\begin{equation}\label{qss}
u_i=\frac{b_ix_i}{\ds b_i+\sum_{j=1}^{n}a_{ij}x_j},
\quad i=1\cdots n \ .
\end{equation}
Since one has $\sum_{j=1}^{n}a_{ij}x_j>0$, for $i=1\cdots n$,
this slow manifold is globally asymptotically stable for
(\ref{fe3}). Thus, the Tikhonov's Theorem  \cite{LST,T,W}
applies and asserts that,
after a fast transition toward the slow manifold, the solutions are approximated by a solution of the reduced equation,
which is obtained by replacing the fast variables $u_i$ in (\ref{eqent3}) by the quasi steady states (\ref{qss}).
One obtains the following reduced model
\begin{equation}\label{reducedN}
\left\{
\begin{array}{l}
\ds \dot{S}=D(S_{in}-S)-\sum_{i=1}^{n}{\mu}_i(S,x)x_i\\
\dot{x}_i=\left({\mu}_i(S,x)-d_{i}(x)\right)x_i, \quad i=1\cdots n
\end{array}
\right.
\end{equation}
 where
\begin{equation}\label{reducedN1}
{\mu}_i(S,x)=f_{i}(S)p_i(x)+g_{i}(S)(1-p_i(x)),\quad
{d}_i(x)=D_{0i}p_i(x)+D_{1i}(1-p_i(x)),\qquad i=1,\cdots,n
\end{equation}
with
\begin{equation}\label{reducedN2}
p_i(x)=\frac{b_i}{\ds b_i+\sum_{j=1}^{n}a_{ij}x_{j}},\qquad x=(x_1,\cdots,x_n).
\end{equation}
Hence, we have shown the following result
\begin{theorem}
Let $(S(t,\varepsilon),u_1(t,\varepsilon),v_1(t,\varepsilon),\cdots,u_n(t,\varepsilon),v_n(t,\varepsilon))$ be the solution of (\ref{Nspecies}) with initial conditions
$S(0)\geqslant 0$ and $u_i(0)>0$, $v_i(0)\geqslant 0$, for $1\leqslant i\leqslant n$.
Let $\left(\overline{S}(t),\overline{x}_1(t),\cdots,\overline{x}_n(t)\right)$ be the solution of the reduced problem (\ref{reducedN})
with initial conditions
$$\overline{S}(0)=S(0),\qquad \overline{x}_i(t)=u_i(0)+v_i(0),\quad 1\leqslant i\leqslant n.$$
Then as $\varepsilon\to 0$
$$S(t,\varepsilon)=\overline{S}(t)+o(1),\qquad x_i(t,\varepsilon)=\overline{x}_i(t)+o(1),\quad 1\leqslant i\leqslant n$$
uniformly for $t\in[0,T]$, and for all $1\leqslant i\leqslant n$, as $\varepsilon\to 0$
$$
u_i(t,\varepsilon)=
\frac{b_i\overline{x}_i(t)}
{b_i+\sum_{j=1}^na_{ij}\overline{x}_j(t)}+o(1)
,\quad
v_i(t,\varepsilon)=
\frac{\left(\sum_{j=1}^na_{ij}\overline{x}_j(t)\right)\overline{x}_i(t)}
{b_i+\sum_{j=1}^na_{ij}\overline{x}_j(t)}+o(1)
$$
uniformly for $t\in[t_0,T]$, where $T>t_0>0$ are arbitrarily fixed. If the solution of the reduced problem tends to an asymptotically stable equilibrium, then we can put $T=+\infty$ in the approximations given above.
\end{theorem}

Since the planktonic bacteria have a better access to the substrate than the bacteria in flocks one assumes $f_{i}(S)>g_{i}(S)$.
Notice that one has $\frac{\partial p_i}{\partial x_j}<0$ for any $i$, $j$. Hence
$$\frac{\partial{\mu}_i}{\partial x_j}=(f_{i}(S)-g_{i}(S))\frac{\partial p_i}{\partial x_j}<0$$
with
$$\frac{\partial{\mu}_i}{\partial S}=f_{i}'(S)p_i(x)+g_{i}'(S)(1-p_i(x))>0 \ .$$
As for the one species case this approach give a motivation to density dependent growth function models, that may lead to species coexistence \cite{LH06,LMR04,LMR06}.

\section{Study of the reduced model for  several species}
\begin{figure}[htbp]
\setlength{\unitlength}{1.0cm}
\begin{center}
\begin{picture}(6.5,4.5)(0,0)
\put(0,-0.5){\rotatebox{0}{\includegraphics[width=8.7cm,height=9cm]{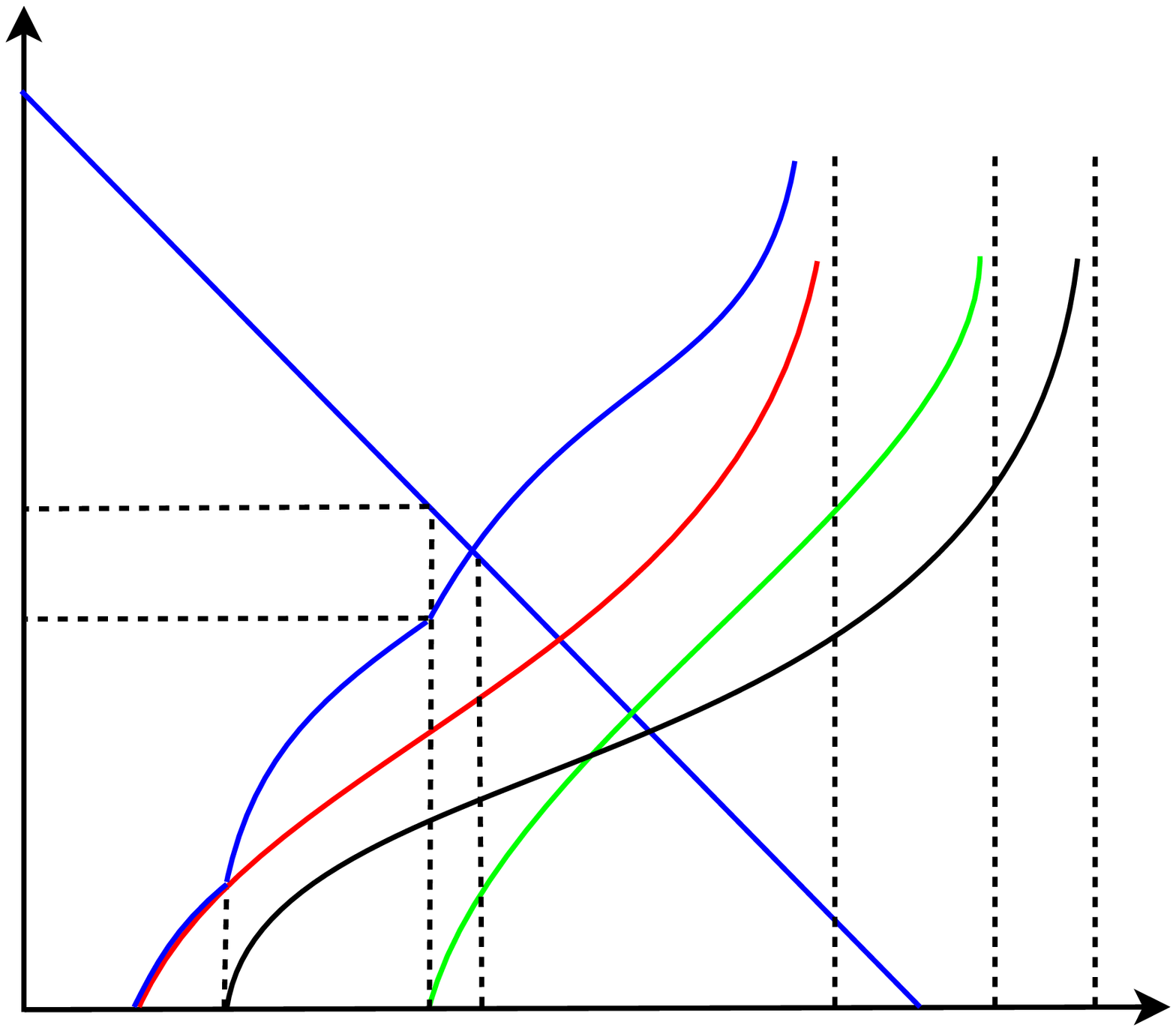}}}
  \put(-0.3,4.3){{\small $DS_{in}$}}
  \put(-1.4,2.3){{\small $D(S_{in}-\lambda_{03})$}}
  \put(-2,1.7){{\small $h_1(\lambda_{03})+h_2(\lambda_{03})$}}
  \put(1.3,3){{\small ${\color{blue}{D(S_{in}-S)}}$}}
  \put(3.3,3.7){{\small ${\color{blue}{\sum_{i=1}^3 h_i(S)}}$}}
  \put(3.8,2.44){{\small ${\color{red}{h_1(S)}}$}}
  \put(5.5,3.6){{\small ${\color{green}{h_3(S)}}$}}
  \put(7.1,3.6){{\small $h_2(S)$}}
  \put(1.1,-0.4){{\small $\lambda_{01}$}}
  \put(1.7,-0.4){{\small $\lambda_{02}$}}
  \put(2.8,-0.4){{\small $\lambda_{03}$}}
  \put(3.3,-0.4){{\small $S^*$}}
  \put(5.2,-0.4){{\small $\lambda_{11}$}}
  \put(5.8,-0.4){{\small $S_{in}$}}
  \put(6.3,-0.4){{\small $\lambda_{13}$}}
  \put(7,-0.4){{\small $\lambda_{12}$}}
  \put(7.6,-0.1){{\small $S$}}
  \put(3.22,2.07){{${\color{blue}{\bullet}}$}}
  \put(3.25,1.37){{${\color{red}{\bullet}}$}}
  \put(3.25,0.87){{$\bullet$}}
  \put(3.25,0.4){{${\color{green}{\bullet}}$}}
\end{picture}
\end{center}
\caption{Condition of existence of the positive equilibrium of (\ref{eq1}) for $n=3$.} \label{fig4}
\end{figure}

The results of this section apply to the reduced model (\ref{reducedN},\ref{reducedN1},\ref{reducedN2}) in the particular case where the rate of attachment and detachment
of species $x_i$ with species $x_j$ are negligible for $i\neq j$, that is to say $a_{ij}=0$ for $i\neq j$. In that case
the function $p_i$ depends only on $x_i$ and is given by
$$
p_i(x)=\frac{b_i}{\ds b_i+a_{ii}x_{i}}
$$
so that the growth function $\mu_i$ and removal rates $d_i$ in (\ref{reducedN}) depend only on $x_i$. We consider then the model
\begin{equation}\label{eq1}
\left\{
\begin{array}{lllll}
\dot{S}  &=& \displaystyle D(S_{in}-S)-\sum_{i=1}^n \mu_i(S,x_i)x_i &&\\[3mm]
\dot{x}_i  &=&  [\mu_i(S,x_i)-d_i(x_i)]x_i &&i=1,\cdots, n
\end{array}
\right.
\end{equation}
This model was studied in \cite{LH06}, in the case when $d_i(x_i)=D$. We do not assume that the functions $\mu_i(S,x_i)$ and $d_i(x_i)$ are of the particular form given by (\ref{reducedN1}). We assume that
\begin{description}
\item[H5:] $\mu_i(0,x_i)=0$ and $\mu_i(S,x_i)\geqslant0$ for all $S > 0$ and all $x_i \geqslant 0.$
\end{description}
\begin{description}
\item[H6:] $\displaystyle\frac{\partial\mu_i}{\partial S}>0$ and $\displaystyle \frac{\partial\mu_i}{\partial x_i}<0$ for all $S>0$ and all $x_i\geqslant 0$.
\end{description}
\begin{description}
\item[H7:] $d_i(0)=D_{0i}$, $d_i(+\infty)=D_{1i}<D_{0i}\leqslant D$, $d_i(x_i)>0$, $d'_i(x_i)<0$ and $[x_id_i(x_i)]'>0$  for all $x_i\geqslant0$.
\end{description}
Let us denote by
$$f_i(S)=\mu_{i}(S,0)\quad\mbox{and}\quad g_i(S)=\mu_{i}(S,+\infty).
$$
The functions $f_i(.)$ and $g_i(.)$ are increasing and positive for all $S > 0$.
If equations $f_{i}(S) = D_{0i}$ and $g_{i}(S) = D_{1i}$
have solutions, one let
$$\lambda_{0i}=f_{i}^{-1}(D_{0i})\quad\mbox{and}\quad \lambda_{1i}=g_{i}^{-1}(D_{1i})$$
otherwise one let $\lambda_{ki}=+\infty$, $k=0,1$.
As for the case of one species (see Assumption {\bf H3}), we add the following assumption
\begin{description}
\item[H8:] $\lambda_{0i}<\lambda_{1i}$ for $i=1\cdots n$. For all $S\in ]\lambda_{0i}, \lambda_{1i}[$ and $x_i\geqslant0$,
one has $d'_i(x_i)>\frac{\partial \mu_i}{\partial x_i}(S,x_i)$.
\end{description}
If the inequality $\lambda_{0i}<\lambda_{1i}$ is reversed for some $i=1\cdots n$, then the situation is much more difficult and will be studied in the future.
Denote
$$ \tilde{\lambda}_{0}=\max \{\lambda_{0i}; \quad i=1,\cdots,n \} \quad\mbox{and}\quad \tilde{\lambda}_1=\min \{\lambda_{1i}; \quad i=1,\cdots,n \}.$$
We assume that
\begin{description}
\item[H9:] $\tilde{\lambda}_{0}< \min(\tilde{\lambda}_1, S_{in})$.
\end{description}
We consider here the existence of a positive equilibrium.
The equilibria of (\ref{eq1}) are solutions of the set of equations
\begin{eqnarray} \label{eq2}
\left\{
\begin{array}{lllll}
\displaystyle D(S_{in}-S)=\sum_{i=1}^n \mu_i(S,x_i)x_i &&\\[3mm]
\mu_i(S,x_i)=d_i(x_i) \quad \mbox{or}\quad x_i=0 &&i=1,\cdots, n.
\end{array}
\right.
\end{eqnarray}
Thus we have to solve the equations
$$\mu_i(S,x_i)=d_i(x_i).$$
Since {\bf H6}, by the implicit function theorem, this equation gives a function  $S=\phi_i(x_i)$
defined for all $x_i\geqslant 0$,
such that
$\phi_i(0)=\lambda_{0i}$, $\phi_i(+\infty)=\lambda_{1i}$ and
$$
\phi'_i(x_i)=\frac{\displaystyle d'_i(x_i)-\frac{\partial \mu_i}{\partial x_i}(S,x_i)}{\displaystyle \frac{\partial \mu_i}{\partial S}(S,x_i)}>0.
$$
The sign of $\phi'(\cdot)$ is given by assumptions {\bf H8}.
We define the function $X_i: S\mapsto X_i(S)$ on $[0,\lambda_{1i}[$ by
$$
\begin{array}{llll}
X_i:[0,&\lambda_{1i}[&\longrightarrow &\mathbb{R}_+\\
                  &S &\longrightarrow & X_i(S)=
                  \left\{
                  \begin{array}{ll}
                    0 & \quad \mbox{if} \quad 0 \leqslant S \leqslant \lambda_{0i}\\
                  x_i=\phi^{-1}_i(S) &\quad \mbox{if} \quad \lambda_{0i}\leqslant S < \lambda_{1i}.
                  \end{array}
                  \right.
\end{array}
$$
Let $h_i(S)=\mu_i(S,X_i(S))X_i(S)$. Since $X_i(\cdot)$ is increasing over $[\lambda_{0i},\lambda_{1i}[$,
so is $h_i(\cdot)$ over this interval. Indeed, one has
$$h'_i(S) = \left(\frac{\partial \mu_i}{\partial S}+\frac{\partial \mu_i}{\partial x_i}X'_i(S)\right)X_i(S)+\mu_i\left(S,X_i(S)\right)X'_i(S).\\
$$
Moreover, for  $S\in]\lambda_{0i},\lambda_{1i}[$, $\mu_i\left(S,X_i(S)\right)=d_i(X_i(S))$ and
$$
X'_i(S)=\frac{\displaystyle\frac{\partial \mu_i}{\partial S}(S,X_i(S))}{\displaystyle d'_i(X_i(S))-\frac{\partial \mu_i}{\partial x_i}(S,X_i(S))}>0.
$$
Then
$$h'_i(S) = \left[d'_i\left(X_i(S)\right)X_i(S)+d_i\left(X_i(S)\right)\right]X'_i(S), \quad \mbox{for } \quad S\in]\lambda_{0i},\lambda_{1i}[ .
$$
Using {{\bf{H7}}},
$$[d_i(x_i)x_i]'=d'_i(x_i)x_i+d_i(x_i)>0,\quad \mbox{for} \quad x_i\geqslant0.
$$
Hence the sign of $h'_i(S)$ is the same as the sign of $X'_i(S)$, that is, $h_i(\cdot)$ is increasing over $[\lambda_{0i},\lambda_{1i}[$
(see Fig. \ref{fig4}).
\\
Consider now the function
$$H(S)=\sum_{i=1}^n h_i(S)-D(S_{in}-S).$$
\begin{lemma}
Equation $H(S)=0$ admits a unique solution $S^*\in]0,\tilde{\lambda}_1[$.
\end{lemma}
{{\bf{Proof}.}}
Since $h_i(S)=0$ for $S\in[0,\lambda_{0i}]$ and $h_i(S)$ is increasing over $[\lambda_{0i},\lambda_{1i}[$ (see Fig. \ref{fig4}),
the function $H(\cdot)$ is increasing over
$(0,\tilde{\lambda}_1)$, and
$$H(0)=-DS_{in}<0 \quad \mbox{and} \quad \lim _{S\rightarrow\tilde{\lambda}_1}H(S)=+\infty.
$$
Hence, there exists a unique $S^*\in]0,\tilde{\lambda}_1[$ such that $H(S^*)=0$.
\medskip
We have the following result :
\begin{proposition}
Assume that {{\bf{H5}}}-{{\bf{H9}}} hold.
System (\ref{eq1}) has a unique positive equilibrium if and only if
\begin{eqnarray}\label{eq3}
\begin{array}{l}
\displaystyle \sum_{i=1}^n \mu_i\Big(\tilde{\lambda}_{0},X_i(\tilde{\lambda}_{0})\Big)X_i(\tilde{\lambda}_{0})<D(S_{in}-\tilde{\lambda}_{0}).
\end{array}
\end{eqnarray}
\end{proposition}
{{\bf{Proof}.}}
A positive equilibrium $E^*=(S^*, x_1^*, \cdots, x_n^*)$, must satisfy
\begin{eqnarray}                                                                                                   \label{eq4}
\begin{array}{l}
\displaystyle D(S_{in}-S^*)=\sum_{i=1}^n \mu_i(S^*,x_i^*)x_i^*
\end{array}
\end{eqnarray}
and
\begin{eqnarray}                                                                                                   \label{eq5}
\begin{array}{l}
\mu_i(S^*,x_i^*)=d_i(x_i^*).
\end{array}
\end{eqnarray}
Equation (\ref{eq5}) is equivalent to $x_i^*=X_i(S^*)$. Thus, (\ref{eq4}) can be written
$$
D(S_{in}-S^*)=\sum_{i=1}^n \mu_i(S^*,X_i(S^*))X_i(S^*)=\sum_{i=1}^n h_i(S^*),
$$
that is  $H(S^*)=0$.
Since $\sum_{i=1}^n h_i(S^*)>0$, then one must have
$$
S^*<S_{in}
\mbox{ and }
S^*>\tilde{\lambda}_{0}.
$$
Notice that $\tilde{\lambda}_{0}<S^*< \tilde{\lambda}_1$ and $S^*<S_{in}$ are satisfied if {{\bf{H9}}} holds.
Then, since $H(S)$ est increasing over $[0,\tilde{\lambda}_1[$,
$$
\tilde{\lambda}_{0}<S^* \Longleftrightarrow H(\tilde{\lambda}_{0})<H(S^*)=0.
$$
Therefore there exists a unique positive equilibrium $S^*$
exactly when $H(\tilde{\lambda}_{0})<0$, which is equivalent to (\ref{eq3}).\\

We study now the asymptotic behavior of the positive equilibrium. First, we establish the following result :
\begin{lemma}\label{PR}
Consider the matrix
\begin{equation}\label{jacobian}
{A}=
\left[
\begin{array}{ccccc}
-D-\sum_{i=1}^n a_i &    c_1 &  c_2   &   \cdots  &  c_n    \\
      a_1           &  -b_1  &   0    &   \cdots  &   0     \\
      a_2           &   0    &  -b_2  &  \cdots   &   0     \\
     \vdots         & \vdots & \vdots &  \ddots   & \vdots  \\
      a_n           &   0    &   0    &  \cdots   &  -b_n
\end{array}
\right]
\end{equation}
Assume that $D>0$ and for $i=1\cdots n$,  $a_i\geqslant 0$, $b_i> 0$ and $c_i\leqslant b_i$.
Then all eigenvalues of $A$ have negative real part.
\end{lemma}

\noindent
{\bf Proof}.
Let $\lambda$ be an eigenvalue of $A$ and $V=(v_0,v_1,\cdots,v_n)\neq 0$ a corresponding eigenvector. Hence we have
\begin{eqnarray} \label{eq29}
\left\{
\begin{array}{rllllll}
\displaystyle \left(-D-\sum_{i=1}^n a_i\right)v_0+\sum_{i=1}^n c_iv_i &=&\lambda v_0 &&\\[5mm]
a_iv_0-b_iv_i &=&\lambda v_i &&i=1,\cdots, n.
\end{array}
\right.
\end{eqnarray}
Assume that $\alpha={\rm Re}(\lambda)\geqslant 0$. Since
$b_i>0$, then $\lambda +b_i\neq 0$.
Therefore,
$$v_i=\frac{a_iv_0}{\lambda+b_i}\quad\mbox{ for }i=1\cdots n.
$$
If $v_0=0$ then $v_i=0$ for $i=1\cdots n$, so that $V=0$ which is impossible.
Thus $v_0\neq 0$ and from the first equation of (\ref{eq29}), we deduce, after simplification by
$v_0$ that
$$
-D-\sum_{i=1}^n a_i+\sum_{i=1}^n a_i\frac{c_i}{\lambda+b_i} =\lambda.
$$
Let $\lambda=\alpha+{\rm i}\beta$. Taking real part of both sides one obtains
$$\alpha=-D-\sum_{i=1}^n a_i+
\sum_{i=1}^na_i\frac{c_i(b_i+\alpha)}{(b_i+\alpha)^2+\beta^2} \cdot
$$
Since $c_i\leqslant b_i$ and $b_i+\alpha>0$ then
$$\frac{c_i(b_i+\alpha)}{(b_i+\alpha)^2+\beta^2}\leqslant
\frac{b_i(b_i+\alpha)}{(b_i+\alpha)^2+\beta^2}\leqslant
\frac{b_i(b_i+\alpha)}{(b_i+\alpha)^2}=
\frac{b_i}{b_i+\alpha}\leqslant 1.$$
Since $a_i\geqslant 0$ then
$$a_i\frac{c_i(b_i+\alpha)}{(b_i+\alpha)^2+\beta^2}\leqslant a_i\quad\mbox{ for }i=1\cdots n.$$
Hence
$$\alpha\leqslant-D-\sum_{i=1}^n a_i+
\sum_{i=1}^na_i<0
$$
which contradicts $\alpha\geqslant 0$.\\

Then, we state the following result :
\begin{proposition}
If $E^*$ exists, then it is locally exponentially stable.
\end{proposition}

\noindent
{\bf Proof}. Since $\mu_i(S^*,x^*_i)=d_i(x^*_i)$,
the Jacobian of the system (\ref{eq1}) at $E^*$ is of the form (\ref{jacobian}) where
$$
a_i=\frac{\partial \mu_i}{\partial S}(S^*,x_i^*)x^*_i
\quad
b_i=-\frac{\partial \mu_i}{\partial x_i}(S^*,x_i^*)x^*_i+x^*_id'_i(x^*_i),
\quad
c_i=-\frac{\partial \mu_i}{\partial x_i}(S^*,x_i^*)x^*_i-d_i(x^*_i).
$$
Since {{\bf{H6}}}, $a_i>0$.
Since {{\bf{H7}}}, $d_i(x^*_i)+x^*_id'_i(x^*_i)>0$, then $-d_i(x^*_i)<x^*_id'_i(x^*_i)$ and hence $c_i<b_i$.
Since {{\bf{H8}}}, $b_i>0$.
The result follows from Lemma \ref{PR}.
\section{Discussion and conclusion}
In this work, we considered a general model of a bio-process with three
compartments involving the substrate, the planktonic and the attached
biomass densities, respectively. Each compartment of the biomass is
characterized by its own specific growth rate and apparent dilution
rate, generalizing previous models of biofilms (with no dilution rate
for the attached bacteria) or models of perfect flocks (with no growth
rate for aggregated individuals). We have analyzed a class of such
models with planktonic and structured biomass, under the assumption
that attachment and detachment processes are fast compared to the
biological scale. Notice that it is only under this assumption that
the main results of the paper are valid, notably the fact that the
three order model (\ref{genmod}) can be reduced to the second order model (\ref{reduced-general}). If it is not the case, the analysis of the original three order model must be done to establish its qualitative behavior.\\
Our study reveals two main characteristics of this model:
 \begin{itemize}
 \item[1.] the reduced dynamics may exhibit a bi-stable behavior even though each growth function is monotonic. This phenomenon is new and is usually met in the chemostat but when the growth function is non-monotonic (such as the Haldane law);
 \item[2.] for bio-processes in which part of the biomass is under a structured form (in flocks or biofilm), the macroscopic models (with reduced dynamics involving only the aggregated biomass and the substrate) should include a growth rate and an apparent dilution rate that are both density dependent. This result contributes to the actual debate in biotechnological engineering involving bio-processes with structured biomass, where it was not clear whether it is better to modify the growth rate functions or the hydrodynamical terms in the macroscopic equations of the system, to cope with the specificity due to the attachment process.
 \end{itemize}

\bigskip

\noindent {\bf Acknowledgments.}
This work has been supported by the DISCO (Multi-scale modelling bioDIversity Structure COupling in biofilms) project,
granted by the French National Research Agency ANR (AAP215-SYSCOMM-2009).

\end{document}